# SHARP ASYMPTOTICS OF THE FUNCTIONAL QUANTIZATION PROBLEM FOR GAUSSIAN PROCESSES

By Harald Luschgy and Gilles Pagès

*Universität Trier and Université Paris 6*

The sharp asymptotics for the $L^2$-quantization errors of Gaussian measures on a Hilbert space and, in particular, for Gaussian processes is derived. The condition imposed is regular variation of the eigenvalues.

**1. Introduction.** The quantization of probability distributions is an old story which starts in the late 1940s. It has been conceived in order to drastically cut down the storage of signal data to be analyzed. For a comprehensive survey of the theory of quantization, including its historical development, we refer to Gray and Neuhoff (1998). For the mathematical aspects of quantization, one may consult Graf and Luschgy (2000), and for more applied aspects in the field of information theory and signal processing, the book of Gersho and Gray (1992) is appropriate.

However, only recently rigorous extensions to the functional quantization of continuous-time stochastic processes have been obtained for the Gaussian case. See Luschgy and Pagès (2002), Dereich, Fehringer, Matoussi and Scheutzow (2003), Dereich (2003) and Graf, Luschgy and Pagès (2003). In particular, the order of convergence to zero of the quantization error has been investigated. The main result of this paper is the sharp asymptotics of the $L^2$-quantization error for a large class of Gaussian processes in a Hilbert space framework. This makes the high-resolution theory in this setting as precise as in the finite-dimensional theory.

The framework can be stated as follows. Let $X$ be a centered Gaussian random vector defined on a probability space $(\Omega, \mathcal{A}, \mathbb{P})$ taking its values in a real separable Hilbert space $H$ with scalar product $\langle \cdot, \cdot \rangle$ and norm $\| \cdot \|$. The distribution $\mathbb{P}_X$ of $X$ will be denoted by $P$ to simplify notations. For $n \in \mathbb{N}$, the $L^2$-quantization problem at level $n$ consists in minimizing

$$E \min_{a \in \alpha} \|X - a\|^2$$

over all sets $\alpha \subset H$ with $|\alpha| \le n$, where $|\cdot|$ is for cardinality. The minimal $n$th quantization error of $P$ is then defined by

$$(1.1) \quad e_n = e_n(P) = \inf \left\{ \left( E \min_{a \in \alpha} \|X - a\|^2 \right)^{1/2} : \alpha \subset H, 1 \le |\alpha| \le n \right\}.$$

The $L^2$-error is the most common measure of the performance of a quantization or lossy data compression system mainly for its simplicity.

Let $\alpha \subseteq H$ be a codebook with $|\alpha| \le n$. One easily shows that the best approximation of $X$ by an $\alpha$-valued random vector is achieved by applying the rule of the nearest neighbor which corresponds to the geometric object called Voronoi partition. So, if

$$f = \sum_{a \in \alpha} a \mathbb{1}_{A_a},$$

where $\{A_a : a \in \alpha\}$ is a Borel measurable partition of $H$ such that, for every $a \in \alpha$, $A_a$ is contained in the (closed and convex) Voronoi region

$$\left\{ x \in H : \|x - a\| = \min_{b \in \alpha} \|x - b\| \right\},$$

then

$$E\|X - f(X)\|^2 = E \min_{a \in \alpha} \|X - a\|^2.$$

Thus one arrives at the representation

$$(1.2) \qquad e_n = \inf_f (E\|X - f(X)\|^2)^{1/2},$$

where the infimum is taken over all $n$-quantizing rules $f$, that is, Borel measurable maps $f : H \to H$ with $|f(H)| \le n$.

We address the issue of high-resolution quantization, that is, the performance of $n$-quantizers and the behavior of $e_n$ as $n \to \infty$. Denote by $K_P \subset H$ the reproducing kernel Hilbert space (Cameron–Martin space) associated to $P$. Observe that $\operatorname{supp}(P)$ coincides with the closure of $K_P$. Let $\lambda_1 \ge \lambda_2 \ge \cdots > 0$ be the ordered nonzero eigenvalues of the covariance operator of $P$ (each written as many times as is its multiplicity) and let $\{u_j : j \ge 1\}$ be a corresponding orthonormal basis of $\operatorname{supp}(P)$ consisting of eigenvectors. If $d := \dim K_P < \infty$, then $e_n(P) = e_n(\bigotimes_{j=1}^d \mathcal{N}(0, \lambda_j))$, the minimal $n$th quantization error of $\bigotimes_{j=1}^d \mathcal{N}(0, \lambda_j)$ with respect to the $l_2$-norm on $\mathbb{R}^d$, and thus we can read off the asymptotic behavior of $e_n$ from the high-resolution formula

$$(1.3) \quad \lim_{n\to\infty} n^{1/d} e_n\left(\bigotimes_{j=1}^d \mathcal{N}(0,\lambda_j)\right) = q(d)\sqrt{2\pi}\left(\prod_{j=1}^d \lambda_j\right)^{1/2d}\left(\frac{d+2}{d}\right)^{(d+2)/4},$$

where $q(d)$ is a constant in $(0,\infty)$ depending only on the dimension $d$ [Zador, see Graf and Luschgy (2000)]. Except in dimension $d=1$ and $d=2$, the true value of $q(d)$ is unknown.

Now assume $\dim K_P = \infty$. Consider the Karhunen–Loève expansion

$$(1.4) \qquad\qquad X = \sum_{j=1}^\infty \lambda_j^{1/2} Z_j u_j,$$

where $Z_j = \langle X, u_j\rangle/\lambda_j^{1/2}$, $j \geq 1$, are i.i.d. $\mathcal{N}(0,1)$-distributed random variables. [It is known that the Karhunen–Loève basis is optimal for quantization of Gaussian measures, see Luschgy and Pagès (2002).] Optimal quantization of $X$ at level $n$ consists in approximating it by a certain finite number $d = d(n)$ of coefficients and the $n$-quantization of these coefficients. More precisely, let $g: \mathbb{R}^d \to \mathbb{R}^d$ denote an $n$-optimal quantizer for $\bigotimes_{j=1}^d \mathcal{N}(0,\lambda_j)$ and set

$$f_n(X) := \sum_{j=1}^d (g(\lambda_1^{1/2} Z_1, \ldots, \lambda_d^{1/2} Z_d))_j u_j.$$

Then

$$e_n^2 = E\|X - f_n(X)\|^2 = \sum_{j\geq d(n)+1} \lambda_j + e_n\left(\bigotimes_{j=1}^{d(n)} \mathcal{N}(0,\lambda_j)\right)^2$$

[see Luschgy and Pagès (2002)]. The critical dimension $d(n)$ is small when compared with $n$ but otherwise unknown for $n \geq 3$. [Since $d(n) \leq n-1$, we have $d(1) = 0$ and $d(2) = 1$.]

In this paper we improve some of the results in Luschgy and Pagès (2002) and derive the sharp asymptotics of $e_n$ as $n \to \infty$ analogously to the finite-dimensional case (1.3) and with slower rates than any $n^{-a}$, of course (Theorems 2.1 and 2.2). This is achieved for regularly varying eigenvalues. The result obtained is even better than (1.3) since limiting constants can be evaluated.

A simple way of obtaining compression is product quantization. Here the Karhunen–Loève coefficients are individually quantized. Thus considering

$$(1.5) \qquad\qquad f_n^{(1)}(X) := \sum_{j=1}^{m(n)} \lambda_j^{1/2} g_j(Z_j) u_j,$$

where $m = m(n) \in \mathbb{N}$ is suitably chosen and $g_j : \mathbb{R} \to \mathbb{R}$ are $n_j$-optimal quantizers for $\mathcal{N}(0,1)$ with optimally allocated $n_j \in \mathbb{N}$ such that $\Pi_{j=1}^m n_j \leq n$, we further show that, for regularly varying eigenvalues with index $-1$ (the largest possible index i.e. the slowest possible decrease), $f_n^{(1)}$ is asymptotically optimal. This means that

$$\lim_n \frac{e_n^2}{E\|X - f_n^{(1)}(X)\|^2} = 1.$$

Furthermore, one shows that $(E\|X - f_n^{(1)}(X)\|^2)^{1/2}$ does follow a sharp rate of convergence as $n \to \infty$ which in turn is that of $e_n$ (see Theorem 2.1). When the eigenvalues are regularly varying with index $-b < -1$, it turns out that some product quantizers $f_n^{(d)}$ with a similar structure but based on quantizing $d$-dimensional marginal blocks are asymptotically almost optimal for some large values of $d$ and provide sharp asymptotics. Furthermore, it is to be noticed that, in that case, the above scalar product quantizers $f_n^{(1)}$ still achieve the sharp rate of convergence. The induced loss is (asymptotically) within a sometimes small constant multiple of the minimal quantization error. These results are stated in Theorem 2.2.

A famous notion of information theory is Shannon's (1949) $\varepsilon$-entropy (rate-distortion function) of $P$. For $\varepsilon > 0$, it is defined by

$$
\begin{aligned}
R(\varepsilon) &= R_P(\varepsilon) \\
&= \inf\Big\{ H(Q|P \otimes Q_2) : Q \text{ probability} \\
&\qquad\qquad \text{on } H \times H \text{ with first marginal } Q_1 = P \\
&\qquad\qquad \text{and } \int_{H \times H} \|x - y\|^2 \, dQ(x, y) \leq \varepsilon^2 \Big\},
\end{aligned}
\tag{1.6}
$$

where $H(Q|P \otimes Q_2)$ denotes the relative entropy (mutual information)

$$H(Q|P \otimes Q_2) = \int_H \log\left(\frac{dQ}{dP \otimes Q_2}\right) dQ$$

if $Q$ is absolutely continuous with respect to the product of the marginals $P \otimes Q_2$ and equals to $\infty$ otherwise. The simple converse part of the source coding theorem [cf. Berger (1971), Theorem 3.2.2, and Graf and Luschgy (2000), page 163] says that the minimal number $N(\varepsilon)$ of codewords needed in a codebook $\alpha$ such that $E\min_{a \in \alpha} \|X - a\|^2 \leq \varepsilon^2$ satisfies

$$\log N(\varepsilon) \geq R(\varepsilon).$$

[In particular, note that $R(e_n) \leq \log n$.] As an application we obtain that $\log N(\varepsilon)$ is precisely $R(\varepsilon)$ in the small distortion regime, that is,

$$\lim_{\varepsilon \to 0} \frac{\log N(\varepsilon)}{R(\varepsilon)} = 1$$

(Corollary 2.4). This sharp asymptotics of the rate of $\log N(\varepsilon)$ is also touched by Donoho (2000).

A further application concerns the small ball problem and its relation to Shannon's $\varepsilon$-entropy.

The paper is organized as follows. In the next section we state the results outlined above. Section 3 contains a collection of examples. Section 4 is devoted to the proofs.

Throughout, all logarithms are natural logarithms and $[x]$ denotes the integer part of the real number $x$.

**2. Statement of results.** Now we formulate sharp asymptotic results for the $n$th quantization errors $e_n = e_n(P)$ and determine the asymptotic behavior of the optimal product quantizers $f_n^{(1)}$ as $n \to \infty$ for centered Gaussian measures $P$ with $\dim K_P = \infty$. It is convenient to use the symbols $\sim$ and $\lesssim$, where $a_n \sim b_n$ means $a_n/b_n \to 1$ and $a_n \lesssim b_n$ means $\limsup_n a_n/b_n \leq 1$.

Let us first give a precise definition of $f_n^{(1)}$. Given $n, m \in \mathbb{N}$, let $n_1, \ldots, n_m \in \mathbb{N}$ with $\Pi_{j=1}^m n_j \leq n$ and let $g_j : \mathbb{R} \to \mathbb{R}$ be $n_j$-optimal quantizers for $\mathcal{N}(0,1)$, $j \in \{1, \ldots, m\}$. Set

$$f(x) = \sum_{j=1}^m \lambda_j^{1/2} g_j(x_j) u_j, \qquad x \in H,$$

where $x_j = \langle x, u_j \rangle / \lambda_j^{1/2}$. Then $|f(H)| \leq n$ and, for every $m \in \mathbb{N}$,

$$E\|X - f(X)\|^2 = \sum_{j \geq m+1} \lambda_j + \sum_{j=1}^m \lambda_j E(Z_j - g_j(Z_j))^2$$

$$= \sum_{j \geq m+1} \lambda_j + \sum_{j=1}^m \lambda_j e_{n_j}(\mathcal{N}(0,1))^2$$

$$\leq \sum_{j \geq m+1} \lambda_j + C(1) \sum_{j=1}^m \lambda_j n_j^{-2},$$

where

$$C(1) := \sup_{k \geq 1} k^2 e_k(\mathcal{N}(0,1))^2$$

is a universal constant.

By the Zador theorem [cf. (1.3)], $C(1) < \infty$. Finally,

$$E\|X - f(X)\|^2 \leq \inf_m \inf_{n_1 \times \cdots \times n_m \leq n} \left( \sum_{j \geq m+1} \lambda_j + C(1) \sum_{j=1}^m \lambda_j n_j^{-2} \right).$$

We may first optimize the integer bit allocation given by the $n'_j s$ for a given $m$ and then select some $m = m(n)$ (hopefully close to the optimal one). To this end, first note that, for a fixed $m \in \mathbb{N}$, the continuous bit allocation problem reads, for every $n \in \mathbb{N}$,

$$\inf\left\{\sum_{j=1}^m \lambda_j y_j^{-2} : y_j > 0, \prod_{j=1}^m y_j \le n\right\} = \sum_{j=1}^m \lambda_j z_j^{-2} = n^{-2/m} m \left(\prod_{j=1}^m \lambda_j\right)^{1/m},$$

with

$$z_j = n^{1/m} \lambda_j^{1/2} \left(\prod_{j=1}^m \lambda_j\right)^{-1/2m}, \qquad j = 1, \ldots, m.$$

One can produce an integer-valued (approximate) solution by setting $n_j = [z_j]$ provided all $z_j \ge 1$. In fact, since the sequence $(\lambda_j)$ is nonincreasing $z_1 \ge \cdots \ge z_m$, so one simply needs that $z_m = n^{1/m} \lambda_m^{1/2} (\Pi_{j=m}^m \lambda_j)^{-1/2m} \ge 1$. A natural choice for the dimension $m$ is then

$$(2.1) \qquad m = m(n) := \max\left\{k \ge 1 : n^{1/k} \lambda_k^{1/2} \left(\prod_{j=1}^k \lambda_j\right)^{-1/2k} \ge 1\right\},$$

$$(2.2) \qquad n_j = n_j(n) := \left[n^{1/m} \lambda_j^{1/2} \left(\prod_{i=1}^m \lambda_j\right)^{-1/2m}\right], \qquad j \in \{1, \ldots, m\},$$

$$(2.3) \quad f_n^{(1)}(x) := \sum_{j=1}^{m(n)} \lambda_j^{1/2} g_j(x_j) u_j, \qquad\qquad\qquad x \in H.$$

We need the notion of a regularly varying function. A measurable function $\varphi : (s, \infty) \to (0, \infty)$ $(s \ge 0)$ is said to be regularly varying at infinity with index $b \in \mathbb{R}$ if, for every $t > 0$,

$$\lim_{x \to \infty} \frac{\varphi(tx)}{\varphi(x)} = t^b.$$

Regular variation of $\varphi : (0, s) \to (0, \infty)$ $(s > 0)$ at zero is defined analogously. Slow variation corresponds to $b = 0$.

THEOREM 2.1. *Assume* $\lambda_j \sim \varphi(j)$ *as* $j \to \infty$, *where* $\varphi : (s, \infty) \to (0, \infty)$ *is a decreasing, regularly varying function at infinity of index* $-1$ *for some* $s \ge 0$. *Set, for every* $x > s$,

$$\psi(x) := \frac{1}{\int_x^\infty \varphi(y)\, dy}.$$

*Then*

$$e_n \sim (E\|X - f_n^{(1)}(X)\|^2)^{1/2} \sim \psi(\log n)^{-1/2} \qquad as\ n \to \infty.$$

*Moreover,*

$$m(n) \sim 2 \log n.$$

REMARK 2.1. Since $\sum_{j=1}^{\infty} \lambda_j < \infty$, the integral $\int_x^{\infty} \varphi(y) \, dy$ is finite. Observe also that the above function $\psi$ is slowly varying at infinity [see Bingham, Goldie and Teugels (1987), Proposition 1.5.9 b]. The most prevalent form for $\varphi$ is

$$\varphi(x) = cx^{-1}(\log x)^{-a}, \qquad a > 1, c > 0, x > 1.$$

Then $\psi(x) = (a-1)(\log x)^{a-1}/c$ and hence

$$e_n \sim (E\|X - f_n^{(1)}(X)\|^2)^{1/2} \sim \left(\frac{c}{a-1}\right)^{1/2} (\log \log n)^{-(a-1)/2} \qquad \text{as } n \to \infty.$$

The following theorem is devoted to the case of regularly varying eigenvalues with index $-b < -1$. It includes a wide class of Gaussian processes. As mentioned in the Introduction, the sharp asymptotics for $e_n$ in item (a) is now approximately achieved by $d$-dimensional marginal block product quantizers. They will be more precisely defined further on in the proof [see (4.1)–(4.3) in Section 4]. Asymptotic optimality is obtained for some high values $d$ of the marginal block dimension. Furthermore, we show in item (b) that, although no longer asymptotically optimal like in Theorem 2.1 (setting $-b = -1$), the scalar product quantizers $f_n^{(1)}$ as defined in (1.5) still provide the sharp rate of convergence to 0 for $e_n$.

THEOREM 2.2. *Assume* $\lambda_j \sim \varphi(j)$ *as* $j \to \infty$, *where* $\varphi : (s, \infty) \to (0, \infty)$ *is a decreasing, regularly varying function at infinity of index* $-b < -1$ *for some* $s \geq 0$. *Set, for every* $x > s$,

$$\psi(x) := \frac{1}{x\varphi(x)}.$$

(a) *Sharp asymptotics for* $e_n$. *Then*

$$e_n \sim \left(\left(\frac{b}{2}\right)^{b-1} \frac{b}{b-1}\right)^{1/2} \psi(\log n)^{-1/2} \qquad \text{as } n \to \infty.$$

(b) *Asymptotics of the scalar product quantizers* $f_n^{(1)}$. *Moreover,*

$$m(n) \sim \frac{2 \log n}{b}$$

*and*

$$\begin{aligned}
(E\|X &- f_n^{(1)}(X)\|^2)^{1/2} \\
&\lesssim \left(\left(\frac{b}{2}\right)^{b-1} \left(\frac{1}{b-1} + 4C(1)\right)\right)^{1/2} \psi(\log n)^{-1/2} \qquad \text{as } n \to \infty,
\end{aligned}$$

*where the real constant $C(1)$ is given by* (2.1).

The proof combines finite-dimensional quantization theory and Shannon's rate-distortion theory.

REMARK 2.2.   (i) We obtain from Theorem 2.2

$$\frac{(E\|X - f_n^{(1)}(X)\|^2)^{1/2}}{e_n} \lesssim \left(\frac{1 + 4C(1)(b-1)}{b}\right)^{1/2} \qquad \text{as } n \to \infty.$$

So, if the index $b$ is close to 1, $f_n^{(1)}$ is close to asymptotic optimality. The constant $C(1)$ is lower bounded by

$$\lim_{k \to \infty} k^2 e_k(\mathcal{N}(0,1))^2 = \frac{\sqrt{3}\pi}{2} = 2.7206\dots$$

[see Graf and Luschgy (2000), page 124]. There is strong numerical evidence for $C(1) = \frac{\sqrt{3}\pi}{2}$. We computed upper bounds of $e_k(\mathcal{N}(0,1))^2$ using the $\frac{i}{k+1}$-quantiles, $i \in \{1, \dots, k\}$, of $\mathcal{N}(0,3)$ which are known to be asymptotically optimal. We thus found

$$\sup_{1 \le k \le 1000} k^2 e_k(\mathcal{N}(0,1))^2 \le \frac{\sqrt{3}\pi}{2}.$$

This suggests that the product quantizing rule $f_n^{(1)}$ cannot be dramatically improved upon for regularly varying eigenvalues.

(ii) The most prevalent form for $\varphi$ is

$$\varphi(x) = c x^{-b} (\log x)^{-a}, \qquad b > 1, a \in \mathbb{R}, x > \max\{1, e^{-a/b}\}.$$

Then we have from the above that

$$e_n \sim \left(c \left(\frac{b}{2}\right)^{b-1} \frac{b}{b-1}\right)^{1/2} (\log n)^{-(b-1)/2} (\log \log n)^{-a/2} \qquad \text{as } n \to \infty.$$

A useful equivalence principle can be deduced from the preceding theorems.

COROLLARY 2.3.   *Assume the situation of Theorem* 2.1 *or* 2.2. *Let $V$ and $W$ be centered Gaussian measures on $H$ and assume that* $\dim \operatorname{supp}(V) < \infty$ *and $W$ is equivalent to $P * V$. Then*

$$e_n(W) \sim e_n(P) \qquad \text{as } n \to \infty.$$

Now we consider, for $\varepsilon > 0$,

$$(2.4) \qquad N(\varepsilon) := \min\{n \geq 1 : e_n \leq \varepsilon\},$$

and the announced strong equivalence of $\log N(\varepsilon)$ and $R(\varepsilon)$. The following "flooding" formula for the $\varepsilon$-entropy $R(\varepsilon)$ of Gaussian measures was originally given by Kolmogorov (1956) [see also Ihara (1993), Theorem 6.9.1]. For $0 < \varepsilon < e_1 = (\sum_{j=1}^{\infty} \lambda_j)^{1/2}$, let

$$(2.5) \qquad r = r(\varepsilon) := \max\left\{ k \geq 1 : \sum_{j \geq k+1} \lambda_j + k\lambda_k > \varepsilon^2 \right\}$$

and let $\vartheta = \vartheta(\varepsilon) \in [\lambda_{r+1}, \lambda_r)$ be uniquely defined by

$$(2.6) \qquad \sum_{j \geq r+1} \lambda_j + r\vartheta = \varepsilon^2.$$

Then $\sum_{j=1}^{\infty} \min\{\lambda_j, \vartheta\} = \varepsilon^2$ and

$$(2.7) \qquad R(\varepsilon) = \tfrac{1}{2} \sum_{j=1}^{r} \log(\lambda_j/\vartheta).$$

The "reproducing distribution" $Q_2 = Q_2(\varepsilon) = \mathbb{P}^Y$ is given by

$$(2.8) \qquad Y = Y(\varepsilon) := \sum_{j=1}^{r} \left[ \lambda_j^{1/2} \left( 1 - \frac{\vartheta}{\lambda_j} \right) Z_j + \vartheta^{1/2} \left( 1 - \frac{\vartheta}{\lambda_j} \right)^{1/2} Z_j' \right] u_j,$$

where $Z_1', Z_2', \ldots$ are i.i.d. $\mathcal{N}(0,1)$-distributed random variables independent of $(Z_j)_{j \geq 1}$, since $Q = Q(\varepsilon) := \mathbb{P}^{(X,Y)}$ solves the minimum problem in Shannon's $R(\varepsilon)$, that is, $R(\varepsilon) = H(Q | P \otimes Q_2)$.

COROLLARY 2.4. *Assume the situation of Theorem* 2.2. *Then*

$$\log N(\varepsilon) \sim R(\varepsilon) \qquad as \ \varepsilon \to 0,$$

$$R(e_n) \sim \log n,$$

$$r(e_n) \sim \frac{2 \log n}{b} \qquad as \ n \to \infty.$$

*Furthermore, $R$ is regularly varying at zero with index $-2/(b-1)$.*

REMARK 2.3. (i) Donoho (2000) states $\log N(\varepsilon) \sim R(\varepsilon)$ for eigenvalues $\lambda_j \sim j^{-b}$ with $b > 1$ and argues that this sharp asymptotics is a consequence of Shannon's rate-distortion theory. Our proof of (the more general) Corollary 2.4 is not in the range of the Shannon theory (see Remark 4.1) and therefore, it does not support Donoho's assessment.

(ii) Let $f_n$ be an $n$-optimal quantizer for $P$. Then, under the condition of Theorem 2.2,

$$\text{entropy}(P^{f_n}) \sim \log n \qquad \text{as } n \to \infty.$$

This follows from the above result, since

$$R(e_n) \leq \text{entropy}(P^{f_n}) \leq \log n.$$

(iii) Since $r(e_n) = \dim \text{supp}(Q_2(e_n))$, the number $r(e_n)$ plays the role of a dimension of the level-$n$ quantization problem. The same role is played by $m(n)$ for the level-$n$ product quantization problem. By Theorem 2.2 and Corollary 2.4, we have $m(n) \sim r(e_n)$.

(iv) In case $\lambda_j \sim cj^{-b}(\log j)^{-a}$ with $c > 0$, $b > 1$ and $a \in \mathbb{R}$, the Shannon $\varepsilon$-entropy can be computed as

$$R(\varepsilon) \sim \frac{b}{2}\left(\frac{cb}{b-1}\left(\frac{b-1}{2}\right)^a\right)^{1/(b-1)} \varepsilon^{-2/(b-1)} \log(1/\varepsilon)^{-a/(b-1)} \qquad \text{as } \varepsilon \to 0$$

[see, e.g., Binia (1974) and Luschgy and Pagès (2002); cf. also (4.14)].

(v) In Shannon information theory is also introduced the distortion-rate function

$$D(R) := \inf\left\{\left(\int \|x-y\|^2 \, dQ(x,y)\right)^{1/2},\right.$$

$$\left. Q \text{ probability on } H \times H, Q_1 = P \text{ and } H(Q|P \otimes Q_2) \leq R\right\},$$

where $H(Q|P \otimes Q_2)$ classically denotes the relative entropy information as defined in the Introduction. One easily checks that it always satisfies $D(\log n) \leq e_n$. Furthermore, under the assumptions of Corollary 2.4, one shows as for the rate-distortion function that

$$D(\log n) \sim e_n \qquad \text{as } n \to \infty.$$

A further application concerns the small ball problem where one tries to find the asymptotic behavior of the function

$$(2.9) \qquad F(\varepsilon) = -\log \mathbb{P}(\|X\| \leq \varepsilon)$$

for small $\varepsilon > 0$. The Shannon $\varepsilon$-entropy provides an upper bound.

COROLLARY 2.5. *Assume the situation of Theorem 2.2. Then*

$$F(\varepsilon) \lesssim R(\varepsilon) \qquad \text{as } \varepsilon \to 0.$$

REMARK 2.4. Under the same condition as above, the lower estimate

$$\left(\frac{b}{b+1}\right)^{b/(b-1)} R(\varepsilon) \lesssim F(\varepsilon) \qquad \text{as } \varepsilon \to 0$$

follows from Theorem 2.5 in Dereich (2003). Simple examples (e.g., Brownian motion and $H = L^2([0,1], dt)$) show that neither $F(\varepsilon) \sim R(\varepsilon)$ nor $F(\varepsilon) \sim (\frac{b}{b+1})^{b/(b-1)} R(\varepsilon)$ as $\varepsilon \to 0$ is true.

**3. Examples.** We consider centered $L^2(\mathbb{P})$-continuous Gaussian processes $X = (X_t)_{t \in I}$ with $I = [0,1]^d$. Then $X$ can be seen as a centered Gaussian random vector with values in the Hilbert space $H = L^2(I, dt)$.

3.1. *Stationary Gaussian processes, Ornstein–Uhlenbeck process and fractional Ornstein–Uhlenbeck processes.* Let $X = (X_t)_{t \in [0,1]}$ be a centered stationary Gaussian process (restricted to $[0,1]$) with covariance function $C(s,t) = \gamma(s-t)$, where $\gamma : \mathbb{R} \to \mathbb{R}$ is continuous, symmetric and positive definite. Assume that the spectral measure admits a (symmetric) Lebesgue density $h$ so that

$$C(s,t) = \int_{\mathbb{R}} e^{i\lambda(t-s)} h(\lambda) \, d\lambda = \int_{\mathbb{R}} \cos(\lambda(t-s)) h(\lambda) \, d\lambda.$$

THEOREM 3.1 [Rosenblatt (1963)]. *Under the condition $h \in L^2(\mathbb{R}, d\lambda)$ (where $d\lambda$ denotes the Lebesgue measure on the real line) and the high-frequency condition*

$$(3.1) \qquad h(\lambda) \sim c\lambda^{-b} \qquad \text{as } \lambda \to +\infty$$

*for some $c > 0$, $b > 1$, the asymptotic behavior of the eigenvalues of the covariance operator is as follows:*

$$(3.2) \qquad \lambda_j \sim 2c\pi^{-(b-1)} j^{-b} \qquad \text{as } j \to \infty.$$

Therefore,

$$(3.3) \qquad e_n \sim \left(2c \left(\frac{b}{2\pi}\right)^{b-1} \frac{b}{b-1}\right)^{1/2} (\log n)^{-(b-1)/2} \qquad \text{as } n \to \infty.$$

Condition (3.1) comprises a broad class of one-dimensional processes including processes with rational spectral densities, the Matérn class [see the discussion in Stein (1999)] and fractional Ornstein–Uhlenbeck processes (but excludes, e.g., bandlimited processes).

The fractional Ornstein–Uhlenbeck process with index $\rho \in (0,2)$ corresponds to

$$C(s,t) = \exp(-a|s-t|^\rho), \qquad a > 0.$$

The spectral measure of this process is a symmetric $\rho$-stable distribution. Its Lebesgue density $h$ is (symmetric) continuous and satisfies

$$h(\lambda) \sim c\lambda^{-(1+\rho)} \qquad \text{as } \lambda \to +\infty$$

with

$$c = \frac{a\Gamma(1+\rho)\sin(\pi\rho/2)}{\pi},$$

where $\Gamma$ denotes the gamma function. Consequently,

$$e_n(FOU) \sim \left(\frac{2a\Gamma(\rho)\sin(\pi\rho/2)(1+\rho)}{\pi}\right)^{1/2}\left(\frac{1+\rho}{2\pi}\right)^{\rho/2}(\log n)^{-\rho/2}$$

(3.4)
$$\text{as } n \to \infty.$$

If $\rho = 1$, one gets the standard stationary Ornstein–Uhlenbeck process on $[0, 1]$. In this case,

$$(3.5) \qquad e_n(OU) \sim \frac{2\sqrt{a}}{\pi}(\log n)^{-1/2} \qquad \text{as } n \to \infty.$$

3.2. *Brownian motion, integrated Brownian motions, Gaussian diffusions and fractional Brownian motions.* (i) *Brownian motion.* Let $B = (B_t)_{t\in[0,1]}$ be a standard Brownian motion. Its covariance operator has eigenvalues

$$\lambda_k = (\pi(k - \tfrac{1}{2}))^{-2}, \qquad k \geq 1.$$

This gives

$$(3.6) \qquad e_n(BM) \sim \frac{\sqrt{2}}{\pi}(\log n)^{-1/2} \qquad \text{as } n \to \infty.$$

(ii) *m-integrated Brownian motion.* For $m \in \mathbb{N}$, let $X = (X_t)_{t\in[0,1]}$ be $m$-times integrated Brownian motion:

$$X_t = \int_0^t \int_0^{s_m} \cdots \int_0^{s_1} B_{s_1}\,ds_1 \cdots ds_m$$

$$= \frac{1}{(m-1)!}\int_0^t (t-s)^{m-1}B_s\,ds.$$

Its covariance function reads

$$C(s,t) = \frac{1}{(m!)^2}\int_0^{s\wedge t}(s-r)^m(t-r)^m\,dr.$$

Ritter ([2000](#), page 79) [see also Freedman ([1999](#)) for $m = 1$ and Gao, Hanning and Torcaso ([2003](#))] has derived the asymptotic behavior of the eigenvalues of the covariance operator:

$$\lambda_k \sim (\pi k)^{-(2m+2)} \qquad \text{as } k \to \infty.$$

Theorem 2.2 then implies that

$$(3.7) \quad e_n(IBM_m) \sim \pi^{-(m+1)}(m+1)^{m+1/2}\left(\frac{2m+2}{2m+1}\right)^{1/2}(\log n)^{-(m+1/2)}$$

The general $m$-times integrated BM as considered by Gao, Hanning and Torcaso (2003) exibits the same asymptotics of the eigenvalues and hence of $e_n$.

(iii) *Gaussian diffusion.* Next, let $X$ be the unique solution of the equation

$$dX_t = A(t)X_t \, dt + dB_t, \qquad X_0 = \xi, \qquad t \in [0,1],$$

where $A \in L^2([0,1], dt)$ and $\xi$ is $\mathcal{N}(0, \sigma^2)$-distributed with $\sigma^2 \geq 0$ and independent of $B$. We find the same asymptotics as for $B$:

$$e_n \sim \frac{\sqrt{2}}{\pi}(\log n)^{-1/2} \qquad \text{as } n \to \infty.$$

This follows from Corollary 2.3. One only has to note that in case $\sigma^2 = 0$ (i.e., $\xi = 0$), the distribution of $X$ is equivalent to the Wiener measure, and in case $\sigma^2 > 0$, it is equivalent to the distribution of $\xi + B$.

(iv) *Fractional Brownian motion.* The fractional Brownian motion (FBM) with Hurst exponent $\beta \in (0,1)$ is a centered continuous Gaussian process on $[0,1]$ having the covariance function

$$C(s,t) = \tfrac{1}{2}(s^{2\beta} + t^{2\beta} - |s-t|^{2\beta}).$$

Using the spectral representation

$$C(s,t) = \Re e \int_{\mathbb{R}} (1 - e^{is\lambda})(1 - e^{-it\lambda})\frac{c}{|\lambda|^{1+2\beta}} \, d\lambda,$$

where

$$c = \frac{1}{4\int_0^\infty (1 - \cos\lambda)\lambda^{-(1+2\beta)} \, d\lambda} = \frac{\Gamma(1+2\beta)\sin(\pi\beta)}{2\pi},$$

one shows the following proposition.

PROPOSITION 3.2. *The ordered eigenvalues of the FBM covariance operator satisfy*

$$\lambda_k \sim 2c\pi^{-2\beta}k^{-(1+2\beta)} \qquad \text{as } k \to \infty.$$

PROOF. A different proof has been found independently by Bronski (2003). For the sake of simplicity we denote by the same letter $C$ a covariance function and its associated kernel operator. The method of proof consists in checking that the eigenvalues of the covariance operator $C$ of the FBM are (strongly) equivalent to those of the stationary covariance kernel

$$C_0(s,t) = \int_{[-1,1]^c} e^{i(s-t)\lambda} \frac{c}{|\lambda|^{1+2\beta}} \, d\lambda.$$

Then, the announced result follows straightforwardly from (3.2) since $\lambda_{0,n} \sim 2c\pi^{-2\beta}n^{-(1+2\beta)}$ as $n \to \infty$. To show this equivalence, we will rely on the following comparison lemma [see, e.g., Rosenblatt (1963)]. $\quad\square$

LEMMA 3.3.    *Let $A_1$, $A_2$ be two symmetric completely continuous transformations on a Hilbert space $H$. Denote the $n$th nonnegative eigenvalue of $A_i$, $i = 1, 2$, and $A_1 + A_2$ by $\lambda_{i,n}^+$ and $\lambda_n^+$, respectively, and the $n$th nonpositive eigenvalue of $A_i$, $i = 1, 2$, and $A_1 + A_2$ by $\lambda_{i,n}^-$ and $\lambda_n^-$. Then, for every $n$, $m \geq 1$,*

$$(3.8) \qquad \lambda_{n+m-1}^+ \leq \lambda_{1,n}^+ + \lambda_{2,m}^+ \quad and \quad \lambda_n^+ \geq \lambda_{1,n+m-1}^+ + \lambda_{2,m}^-.$$

First, we decompose $\widetilde{C} := C - C_0$ as follows:

$$\widetilde{C}(s,t) = C_1(s,t) + C_2(s,t)$$

with

$$C_1(s,t) := \Re e \int_{[-1,1]} (1 - e^{-is\lambda})(1 - e^{it\lambda}) \mu(d\lambda)$$

and

$$C_2(s,t) := \int_{[-1,1]^c} (1 - \cos(t\lambda) - \cos(s\lambda)) \mu(d\lambda),$$

where $\mu(d\lambda) := \frac{c}{|\lambda|^{1+2\beta}} \, d\lambda$.

The operator $C_1$ is nonnegative so $\lambda_{1,n} = \lambda_{1,n}^+$, $n \geq 1$. The range of operator $C_2$ has obviously dimension 2 since it maps $L^2([0,1], dt)$ onto $\langle \mathbf{1}, \gamma \rangle$ with $\gamma(t) := \int_{[-1,1]^c} \cos(t\lambda)\mu(d\lambda)$. Hence $\lambda_{2,m}^\pm = 0$ as soon as $m \geq 3$ and it follows from inequalities (3.8) that, for every $n \geq 3$

$$\lambda_{1,n+2} \leq \widetilde{\lambda}_{1,n}^+ \leq \lambda_{1,n-2} \quad and \quad \widetilde{\lambda}_{1,n}^- = 0.$$

To estimate the eigenvalues $\lambda_{1,n}$ of the operator $C_1$, one first notes that

$$C_1(s,t) = -2c \sum_{k \geq 1} \frac{(-1)^k}{(2k)!} (s^{2k} + t^{2k} - (s-t)^{2k}) \int_0^1 \lambda^{2(k-\beta)-1} \, d\lambda$$

$$= -c \sum_{k \geq 1} \frac{(-1)^k}{(2k)!(k-\beta)} (s^{2k} + t^{2k} - (s-t)^{2k})$$

$$= -c \sum_{k=1}^{n} \frac{(-1)^k}{(2k)!(k-\beta)} (s^{2k} + t^{2k} - (s-t)^{2k})$$

$$\underbrace{\phantom{-c \sum_{k=1}^{n} \frac{(-1)^k}{(2k)!(k-\beta)}}}_{=:C_{1,n}(s,t)}$$

$$+ (-c) \sum_{k \geq n+1} \frac{(-1)^k}{(2k)!(k-\beta)} (s^{2k} + t^{2k} - (s-t)^{2k}).$$

$$\underbrace{\phantom{+ (-c) \sum_{k \geq n+1} \frac{(-1)^k}{(2k)!(k-\beta)}}}_{=:C'_{1,n}(s,t)}$$

One checks that $C_{1,n}$ maps $L^2([0,1], dt)$ into the $(2n-1)$-dimensional space $\mathbb{R}_{2(n-1)}[X]$ of polynomial functions $P$ such that $degree(P) \leq 2(n-1)$. The above lemma implies that

$$\lambda_{1,2n} \leq \lambda_{(1,n),2n}^{+} + \lambda_{(1,n),1}^{\prime +} = \lambda_{(1,n),1}^{\prime +} = \sup_{\|u\|=1} \int_{[0,1]^2} C'_{1,n}(s,t) u(s) u(t) \, ds \, dt.$$

Now, using that $|s^{2k} + t^{2k} - (s-t)^{2k}| \leq 3$ for every $s, t \in [0,1]$, one easily derives that

$$0 \leq \lambda_{1,2n} \leq 3c \sum_{k \geq n+1} \frac{1}{(2k)!(k-\beta)} \leq \frac{3c}{k+1-\beta} \sum_{k \geq 2(n+1)} \frac{1}{k!}$$

$$\leq \frac{3c}{n+1-\beta} \frac{1}{(2n+1)(2n+1)!} \sim \frac{3c}{2n^2} \frac{1}{(2n+1)!}.$$

The sequence $(\lambda_{1,n})_{n \geq 1}$ being nonincreasing, it follows that $\lambda_{1,n} = o(\lambda_{0,n})$ and consequently $\widetilde{\lambda}_n^{+} = o(\lambda_{0,n})$. Finally, one derives the announced conclusion from the equality $C = C_0 + \widetilde{C}$ and inequalities (3.8) of the lemma:

$$\lambda_n = \lambda_n^{+} = \lambda_{0,n} + o(\lambda_{0,n}) \qquad \text{as } n \to \infty.$$

This yields

$$e_n(FBM) \sim \left( \frac{\Gamma(2\beta) \sin(\pi\beta)(1+2\beta)}{\pi} \right)^{1/2} \left( \frac{1+2\beta}{2\pi} \right)^{\beta} (\log n)^{-\beta}$$

$$\text{as } n \to \infty.$$
(3.9)

It is interesting to observe that the quantization error of the fractional Brownian motion exhibits the same asymptotic behavior as that of the fractional Ornstein–Uhlenbeck process with covariance $\exp(-|s-t|^{2\beta}/2)$.

3.3. *Gaussian sheets.* We consider centered ($L^2(\mathbb{P})$-continuous) Gaussian fields $X = (X_t)_{t \in [0,1]^d}$ with covariance function of tensor product form

$$C(s,t) = \prod_{j=1}^{d} C_j(s_j, t_j),$$

where $C_j$ are covariance functions on $[0,1]$. Let $\lambda_1(j) \geq \lambda_2(j) \geq \cdots > 0$ and $\lambda_1 \geq \lambda_2 \geq \cdots 0$ denote the ordered nonzero eigenvalues associated to $C_j$ and $C$ respectively. We rely on the following proposition.

PROPOSITION 3.4 [Papageorgiou and Wasilkowski (1990)]. *If $\lambda_k(j) \sim c_j k^{-b}$ as $k \to \infty$ for every $j \in \{1,\ldots,d\}$, where $c_j > 0, b > 1$, then*

$$\lambda_k \sim \left(\prod_{j=1}^{d} c_j\right)((d-1)!)^{-b} k^{-b}(\log k)^{b(d-1)} \qquad \text{as } k \to \infty.$$

• *Fractional Ornstein–Uhlenbeck sheet.* The fractional Ornstein–Uhlenbeck sheet on $[0,1]^d$ with index $\rho \in (0,2)$ corresponds to the following covariance function

$$C(s,t) = \prod_{j=1}^{d} \exp(-a_j |s_j - t_j|^\rho), \qquad a_j > 0.$$

By Example 3.1 and Proposition 3.4, the eigenvalues of its covariance operator satisfy

$$\lambda_k \sim \left(\prod_{j=1}^{d} a_j\right)\left(\frac{2\Gamma(1+\rho)\sin(\pi\rho/2)}{\pi^{1+\rho}}\right)^d \left(\frac{(\log k)^{d-1}}{(d-1)!k}\right)^{1+\rho} \qquad \text{as } k \to \infty.$$

Therefore,

$$\begin{aligned}
e_n(FOUS) \sim {}& \left(\prod_{j=1}^{d} a_j\right)^{1/2} \left(\frac{2\Gamma(1+\rho)\sin(\pi\rho/2)}{\pi^{1+\rho}}\right)^{d/2} ((d-1)!)^{-(1+\rho)/2} \\
& \times \left(\left(\frac{1+\rho}{2}\right)^\rho \frac{1+\rho}{\rho}\right)^{1/2} (\log n)^{-\rho/2}(\log\log n)^{(1+\rho)(d-1)/2} \\
& \hspace{8cm} \text{as } n \to \infty.
\end{aligned}$$

(3.10)

If $\rho = 1$, one gets the stationary Ornstein–Uhlenbeck sheet on $[0,1]^d$. In this case

$$e_n(OUS) \sim \left(\prod_{j=1}^{d} a_j\right)^{1/2} \frac{2^{(d+1)/2}}{\pi^d (d-1)!}(\log n)^{-1/2}(\log\log n)^{d-1}$$

(3.11)

$$\hspace{8cm} \text{as } n \to \infty.$$

The 2-parameter O.U.-sheet has been successfully used as model for image compression [see, e.g., Rosenfeld and Kak (1976)].

• *Fractional Brownian sheet.* The fractional Brownian sheet with Hurst exponent $\beta \in (0,1)$ is a centered continuous Gaussian field on $[0,1]^d$ having the covariance function

$$C(s,t) = 2^{-d} \prod_{j=1}^{d} (s_j^{2\beta} + t^{2\beta} - |s_j - t_j|^{2\beta}).$$

By Propositions 3.2 and 3.4, the eigenvalues of th FBS covariance operator satisfy

$$\lambda_k \sim \left( \frac{\Gamma(1+2\beta)\sin(\pi\beta)}{\pi^{1+2\beta}} \right)^d ((d-1)!)^{-(1+2\beta)} k^{-(1+2\beta)} (\log k)^{(1+2\beta)(d-1)}$$

$$\text{as } k \to \infty.$$

This yields

$$e_n(FBS) \sim \left( \frac{\Gamma(1+2\beta)\sin(\pi\beta)}{\pi^{1+2\beta}} \right)^{d/2} ((d-1)!)^{-(1+2\beta)/2}$$

(3.12)
$$\times \left( \left( \frac{1+2\beta}{2} \right)^{2\beta} \frac{1+2\beta}{2\beta} \right)^{1/2} (\log n)^{-\beta} (\log\log n)^{(1+2\beta)(d-1)/2}$$

$$\text{as } n \to \infty.$$

If $\beta = \frac{1}{2}$, one gets Brownian sheet where

$$C(s,t) = \prod_{j=1}^{d} (s_j \wedge t_j).$$

In this case

(3.13)   $$e_n(BS) \sim \frac{\sqrt{2}}{\pi^d (d-1)!} (\log n)^{-1/2} (\log\log n)^{d-1} \qquad \text{as } n \to \infty.$$

The same asymptotic behavior of the eigenvalues and subsequently of $e_n$ is obtained for the completely tucked Brownian sheet on $[0,1]^d$, where

$$C(s,t) = \prod_{j=1}^{d} (s_j \wedge t_j - s_j t_j).$$

**4. Proofs of results.** We need an extension of the quantizing rule $f_n^{(1)}$ based now on quantizing blocks of Karhunen–Loève coefficients of fixed block length $d$. Fix $d, n \in \mathbb{N}$ and set $\nu_j := \lambda_{(j-1)d+1}$, $j \geq 1$. Let

(4.1)   $$m = m(n,d) := \max\left\{ k \geq 1 : n^{1/k} \nu_k^{d/2} \left( \prod_{j=1}^{k} \nu_j \right)^{-d/2k} \geq 1 \right\},$$

(4.2)   $$n_j = n_j(n,d) := \left[ n^{1/m} \nu_j^{d/2} \left( \prod_{i=1}^{m} \nu_i \right)^{-d/2m} \right], \qquad j \in \{1, \ldots, m\},$$

and

$$f_n^{(d)}(x) := \sum_{j=1}^{m} \sum_{k=1}^{d} \lambda_{(j-1)d+k}^{1/2} (g_j(x_{(j-1)d+1}, \ldots, x_{jd}))_k u_{(j-1)d+k}, \qquad x \in H,$$

where $g_j : \mathbb{R}^d \to \mathbb{R}^d$ denotes an $n_j$-optimal quantizer for $\mathcal{N}(0, I_d)$. Setting $Z^{(j)} := (Z_{(j-1)d+1}, \ldots, Z_{jd})$, we get

$$(4.3) \qquad f_n^{(d)}(X) = \sum_{j=1}^{m} \sum_{k=1}^{d} \lambda_{(j-1)d+k}^{1/2} (g_j(Z^{(j)}))_k u_{(j-1)d+k},$$

$n_j \geq 1, \Pi_{j=1}^{m} n_j \leq n$ and $|f_n^{(d)}(H)| \leq n$. For $d \geq 2$, the procedure (4.3) is worse than quantizing $d$-blocks of coefficients $\lambda_j^{1/2} Z_j$ but good enough for our purpose. For the evaluation of the error of $f_n^{(d)}$, we need the constant

$$(4.4) \qquad C(d) := \sup_{k \geq 1} k^{2/d} e_k(\mathcal{N}(0, I_d))^2.$$

By the Zador theorem [see (1.3)], $C(d) < \infty$.

LEMMA 4.1.  *For every* $d, n \in \mathbb{N}$,

$$E\|X - f_n^{(d)}(X)\|^2 \leq \sum_{j \geq md+1} \lambda_j + 4^{1/d} C(d) m \nu_m$$

*with* $m = m(n, d)$ *from* (4.1) *and* $\nu_m = \lambda_{(m-1)d+1}$.

PROOF.   We have

$$E\|X - f_n^{(d)}(X)\|^2 = \sum_{j \geq md+1} \lambda_j + \sum_{j=1}^{m} \sum_{k=1}^{d} \lambda_{(j-1)d+k} E(Z_{(j-1)d+k} - (g_j(Z^{(j)}))_k)^2$$

$$\leq \sum_{j \geq md+1} \lambda_j + \sum_{j=1}^{m} \nu_j E\|Z^{(j)} - g_j(Z^{(j)})\|^2$$

$$= \sum_{j \geq md+1} \lambda_j + \sum_{j=1}^{m} \nu_j e_{n_j}(\mathcal{N}(0, I_d))^2.$$

Moreover, by (4.2),

$$\sum_{j=1}^{m} \nu_j e_{n_j}(\mathcal{N}(0, I_d))^2 \leq C(d) \sum_{j=1}^{m} \nu_j n_j^{-2/d}$$

$$= C(d) \sum_{j=1}^{m} \nu_j (n_j + 1)^{-2/d} \left( \frac{n_j + 1}{n_j} \right)^{2/d}$$

$$\leq 4^{1/d} C(d) m n^{-2/dm} \left( \prod_{j=1}^{m} \nu_j \right)^{1/m}$$

$$\leq 4^{1/d} C(d) m \nu_m. \qquad \square$$

A Shannon-type lower bound is as follows.

LEMMA 4.2. *For every* $n \in \mathbb{N}$,

$$e_n^2 \geq \sum_{j \geq m+1} \lambda_j + m \lambda_{m+1}$$

*with* $m = m(n)$ *from* (2.1).

[Note that $m(n) = m(n,1)$.]

PROOF. Setting

$$(4.5) \qquad a_k := \frac{k}{2} \log \left( \left( \prod_{j=1}^{k} \lambda_j \right)^{1/k} \Big/ \lambda_k \right) = \frac{1}{2} \sum_{j=1}^{k} \log(\lambda_j / \lambda_k),$$

we see that

$$(4.6) \qquad m = m(n) = \max\{ k \geq 1 : a_k \leq \log n \}.$$

On the other hand, by (2.8) for $\varepsilon < e_1$,

$$R(\varepsilon) > a_{r(\varepsilon)},$$

so that by the converse source coding theorem, for every $n \geq 2$,

$$\log n \geq R(e_n) > a_{r(e_n)}.$$

Consequently, $r(e_n) \leq m(n)$ and, using (2.5), this yields

$$e_n^2 \geq \sum_{j \geq m+2} \lambda_j + (m+1)\lambda_{m+1} = \sum_{j \geq m+1} \lambda_j + m \lambda_{m+1}.$$

The assertion is also true for $n = 1$ since $m = m(1)$ equals the multiplicity of $\lambda_1$ and

$$e_1^2 = \sum_{j=1}^{\infty} \lambda_j = \sum_{j \geq m+1} \lambda_j + m \lambda_m \geq \sum_{j \geq m+1} \lambda_j + m \lambda_{m+1}. \qquad \square$$

PROOF OF THEOREM 2.1. The subsequent arguments already occur (somewhat hidden) in Luschgy and Pagès (2002). We repeat them for completeness and the reader's convenience. By Lemma 4.1, we have, for every $n \in \mathbb{N}$,

$$e_n^2 \leq E \| X - f_n^{(1)}(X) \|^2 \leq \sum_{j \geq m+1} \lambda_j + 4 C(1) m \lambda_m$$

with $m = m(n)$, where the approximation error is the dominating term. In fact, it follows from the assumption on the eigenvalues that, for $a_k$ defined in (4.5), we have $a_k \sim k/2$ as $k \to \infty$ and hence by (4.6),

$$m(n) \sim 2\log n \qquad \text{as } n \to \infty.$$

This yields

$$\sum_{j \geq m+1} \lambda_j \sim \psi(m)^{-1} \sim \psi(\log n)^{-1} \qquad \text{as } n \to \infty.$$

Moreover,

$$x\varphi(x) = o(\psi(x)^{-1}) \qquad \text{as } x \to \infty$$

[cf. Bingham, Goldie and Teugels (1987), Proposition 1.5.9 b] and thus

$$E\|X - f_n^{(1)}(X)\|^2 \lesssim \psi(\log n)^{-1} \qquad \text{as } n \to \infty.$$

The lower estimate

$$e_n^2 \gtrsim \psi(\log n)^{-1} \qquad \text{as } n \to \infty$$

follows from Lemma 4.2. $\quad\square$

Now we turn to the proof of Theorem 2.2. Let $Q(d)$ denote the quantization coefficient (of order 2) of the $d$-dimensional standard normal distribution $\mathcal{N}(0, I_d)$, that is,

$$(4.7) \qquad Q(d) = \lim_{n \to \infty} n^{2/d} e_n(\mathcal{N}(0, I_d))^2$$

[see (1.3)].

PROPOSITION 4.3.   *The sequence* $(Q(d))_{d \geq 1}$ *satisfies* $\lim_{d \to \infty} Q(d)/d = 1$.

PROOF.   See Graf and Luschgy [(2000), Proposition 9.5]. $\quad\square$

The key property is the following $d$-asymptotics of the constants $C(d)$ defined in (4.4).

PROPOSITION 4.4.   *The sequence* $(C(d))_{d \geq 1}$ *satisfies* $\liminf_{d \to \infty} C(d)/d = 1$.

PROOF.   Since $C(d) \geq Q(d)$ for every $d \in \mathbb{N}$, it follows from Proposition 4.3 that

$$(4.8) \qquad \liminf_{d \to \infty} \frac{C(d)}{d} \geq 1.$$

Again by Proposition 4.3, the converse inequality is true if $C(d) = Q(d)$ for all but finitely many $d$'s. So assume that $C(d) > Q(d)$ for all members in a subsequence of $(C(d), Q(d))_d$. (No special notation for subsequences is used.) Set $e_n(d) := e_n(\mathcal{N}(0, I_d))$ and

$$A(d, k) := \frac{k^{2/d} e_k(d)^2}{d}.$$

If for $d \in \mathbb{N}$, $C(d) > Q(d)$ holds, then choose $\eta > 0$ such that $C(d) > Q(d) + \eta$. By (4.7), there exists $k_0 \in \mathbb{N}$ such that

$$\sup_{k \geq k_0+1} k^{2/d} e_k(d)^2 \leq Q(d) + \eta.$$

Consequently, $C(d) = \sup_{k \leq k_0} k^{2/d} e_k(d)^2$ and hence, there exists $p(d) \in \mathbb{N}$ such that

$$\frac{C(d)}{d} = A(d, p(d)).$$

We claim that for every sequence $(k(d))_d$ in $\mathbb{N}$,

$$\liminf_{d \to \infty} A(d, k(d)) \leq 1. \tag{4.9}$$

The proof of (4.9) which settles the proposition is given by a sequence of steps.

STEP 1.   Assume

$$\liminf_{d \to \infty} \frac{\log k(d)}{d} = 0.$$

By taking a subsequence, we may assume that $\lim_{d \to \infty} \frac{\log k(d)}{d} = 0$. Using the rough upper bound

$$e_{k(d)}(d)^2 \leq e_1(d)^2 = d,$$

one gets

$$A(d, k(d)) \leq k(d)^{2/d}.$$

Consequently,

$$\liminf_{d \to \infty} A(d, k(d)) \leq \lim_{d \to \infty} k(d)^{2/d} = 1.$$

STEP 2.   For $\delta > 0$ and $\varepsilon \in (0, 1)$, consider the special sequence

$$k(d) = [\exp(d(\delta + R_1(\varepsilon)))], \tag{4.10}$$

where $R_1(\varepsilon)$ denotes the $\varepsilon$-entropy of $\mathcal{N}(0, 1)$ given by

$$R_1(\varepsilon) = \log(1/\varepsilon). \tag{4.11}$$

The direct part of Shannon's source coding theorem says that

$$\limsup_{d \to \infty} \frac{e_{k(d)}(d)^2}{d} \leq \varepsilon^2$$

[see Dembo and Zeitouni (1998), Theorem 3.6.2]. A careful reading of the proof shows that their large deviation approach also works for the unbounded (squared) error function in our setting. Since $\lim_{d \to \infty} k(d)^{2/d} = (1/\varepsilon)^2 \exp(2\delta)$, one gets

(4.12)                     $$\limsup_{d \to \infty} A(d, k(d)) \leq \exp(2\delta).$$

STEP 3.   Assume

$$0 < \liminf_{d \to \infty} \frac{\log k(d)}{d} < \infty.$$

By taking a subsequence, we may assume that

$$\lim_{d \to \infty} \frac{\log k(d)}{d} = c \qquad \text{with } 0 < c < \infty.$$

Choose $c_1 \in (0, c)$ and $\delta \in (0, c_1)$. Then for $d$ large, $\frac{\log k(d)}{d} \geq c_1$ and hence

$$k(d) \geq \exp(dc_1) = \exp(d\delta + d(c_1 - \delta)).$$

Set $\varepsilon := \exp(\delta - c_1)$ and

$$q(d) := [\exp(d(\delta + R_1(\varepsilon)))].$$

Then for $d$ large , $k(d) \geq q(d)$ and

$$A(d, k(d)) \leq \frac{k(d)^{2/d}}{d} e_{q(d)}(d)^2 = \frac{k(d)^{2/d}}{q(d)^{2/d}} A(d, q(d)).$$

By Step 2, we have

$$\limsup_{d \to \infty} A(d, q(d)) \leq \exp(2\delta).$$

Since $\lim_{d \to \infty} k(d)^{2/d} = \exp(2c)$ and $\lim_{d \to \infty} q(d)^{2/d} = \exp(2c_1)$, one obtains

$$\limsup_{d \to \infty} A(d, k(d)) \leq \exp(2(c - c_1) + 2\delta).$$

Letting $\delta \to 0$ and then $c_1 \to c$ yields

$$\limsup_{d \to \infty} A(d, k(d)) \leq 1.$$

STEP 4. Assume

$$\lim_{d \to \infty} \frac{\log k(d)}{d} = \infty.$$

Fix $m \in \mathbb{N}$ and proceed by a block-quantizer design consisting of $d$ blocks of length $m$ for quantizing $\mathcal{N}(0, I_{md})$. Set

$$s = s(d) := [k(md)^{1/d}].$$

Then $s^d \le k(md)$ and

$$e_{k(md)}(md)^2 \le e_{s^d}(md)^2 \le d e_s(m)^2.$$

Consequently, for every $d \in \mathbb{N}$,

$$A(md, k(md)) \le k(md)^{2/md} e_s(m)^2 \frac{1}{m} = \left( \frac{k(md)^{1/d}}{s} \right)^{2/m} s^{2/m} e_s(m)^2 \frac{1}{m}.$$

Since $s(d) \to \infty$ as $d \to \infty$, it follows from the Zador theorem that

$$\lim_{d \to \infty} s(d)^{2/m} e_{s(d)}(m)^2 = Q(m)$$

[see (4.7)]. This implies

$$\liminf_{d \to \infty} A(md, k(md)) \le \frac{Q(m)}{m}.$$

Using Proposition 4.3, we deduce that

$$\liminf_{d \to \infty} A(d, k(d)) \le 1. \qquad \square$$

REMARK 4.1. We emphasize that Step 4 of the above proof is not in the range of the Shannon theory.

PROOF OF THEOREM 2.2. By Lemma 4.1, we have, for every $d, n \in \mathbb{N}$,

$$e_n^2 \le E\|X - f_n^{(d)}(X)\|^2 \le \sum_{j \ge md+1} \lambda_j + 4^{1/d} C(d) m \nu_m$$

with $m = m(n, d)$. Setting

$$a_k(d) := \frac{k}{2} \log \left( \left( \prod_{j=1}^{k} \nu_j^d \right)^{1/k} \bigg/ \nu_k^d \right),$$

we see that

$$m(n, d) = \max\{k \ge 1 : a_k(d) \le \log n\}.$$

The assumption on the eigenvalues implies $a_k(d) \sim bdk/2$ as $k \to \infty$ and hence, for every $d \in \mathbb{N}$,

$$m(n,d) \sim \frac{2\log n}{bd} \qquad \text{as } n \to \infty.$$

Consequently,

$$\nu_m \sim \lambda_{md} \sim \left(\frac{2}{b}\right)^{-b} \varphi(\log n),$$

$$m, d\nu_m \sim \left(\frac{2}{b}\right)^{1-b} \psi(\log n)^{-1}$$

and

$$\sum_{j \geq md+1} \lambda_j \sim \frac{md\varphi(md)}{b-1} \sim \frac{1}{b-1}\left(\frac{2}{b}\right)^{1-b} \psi(\log n)^{-1} \qquad \text{as } n \to \infty.$$

We deduce that, for every $d \in \mathbb{N}$,

$$E\|X - f_n^{(d)}(X)\|^2 \lesssim \left(\frac{b}{2}\right)^{b-1}\left(\frac{1}{b-1} + \frac{4^{1/d}C(d)}{d}\right)\psi(\log n)^{-1}$$

(4.13)

$$\text{as } n \to \infty.$$

Note that for $d = 1$, (4.13) gives the desired upper estimate for $f_n^{(1)}$. Now it follows from Proposition 4.4 that

$$e_n^2 \lesssim \left(\frac{b}{2}\right)^{b-1} \frac{b}{b-1} \psi(\log n)^{-1} \qquad \text{as } n \to \infty.$$

The lower estimate

$$e_n^2 \gtrsim \left(\frac{b}{2}\right)^{b-1} \frac{b}{b-1} \psi(\log n)^{-1} \qquad \text{as } n \to \infty$$

is a consequence of Lemma 4.2.  $\square$

Finally, we prove Corollaries 2.3–2.5.

PROOF OF COROLLARY 2.3.   Let $\rho_1 \geq \rho_2 \geq \cdots > 0$ denote the nonzero eigenvalues of the covariance operator of $P*V$ and let $d := \dim \operatorname{supp}(V)$. Then, by the minimax characterization [see (3.8)] of eigenvalues, for every $j \in \mathbb{N}$,

$$\rho_{j+d} \leq \lambda_j \leq \rho_j.$$

Regular variation of the eigenvalues $\lambda_j$ implies $\rho_j \sim \lambda_j$ as $j \to \infty$. Let $\mu_1 \geq \mu_2 \geq \cdots > 0$ denote the nonzero eigenvalues associated to $W$. Theorem 2 in Ihara ([1970](#)) and regular variation of $\rho_j$ imply that $\mu_j \sim \rho_j$ as $j \to \infty$. Thus the assertion follows from Theorems 2.1 and 2.2. $\quad\square$

PROOF OF COROLLARY 2.4. By Theorem 4.12 in Luschgy and Pagès ([2002](#)), we have

$$e_n^2 \sim \left(\frac{b}{2}\right)^{b-1} \frac{b}{b-1} \psi(R(e_n))^{-1} \qquad \text{as } n \to \infty,$$

with the function $\psi$ from Theorem 2.2. Combining this with Theorem 2.2 gives

$$\psi(R(e_n)) \sim \psi(\log n) \qquad \text{as } n \to \infty.$$

There exists a function $\widetilde{\psi}$ which is regularly varying at infinity of index $1/(b-1)$ such that

$$\widetilde{\psi}(\psi(x)) \sim x \qquad \text{as } x \to \infty$$

[cf. Bingham, Goldie and Teugels ([1987](#)), Theorem 1.5.12]. Hence

$$R(e_n) \sim \widetilde{\psi}(\psi(R(e_n))) \sim \widetilde{\psi}(\psi(\log n)) \sim \log n \qquad \text{as } n \to \infty.$$

In particular,

$$\log N(\varepsilon) \sim R(e_{N(\varepsilon)}),$$

$$\log N(\varepsilon) \sim \log(N(\varepsilon) - 1) \sim R(e_{N(\varepsilon)-1}) \qquad \text{as } \varepsilon \to 0.$$

Since $\varepsilon < e_{N(\varepsilon)-1}$ for $\varepsilon \le e_2$ and $e_{N(\varepsilon)} \le \varepsilon$ and thus $R(\varepsilon) \ge R(e_{N(\varepsilon)-1})$ and $R(\varepsilon) \le R(e_{N(\varepsilon)})$, we obtain

$$\log N(\varepsilon) \sim R(\varepsilon) \qquad \text{as } \varepsilon \to 0.$$

Using Theorem 2.2, this implies

$$\varepsilon^2 \sim e_{N(\varepsilon)}^2 \sim \left(\frac{b}{2}\right)^{b-1} \frac{b}{b-1} \psi(R(\varepsilon))^{-1} \qquad \text{as } \varepsilon \to 0.$$

Consequently,

$$R(\varepsilon) \sim \widetilde{\psi}\left(\left(\frac{b}{2}\right)^{b-1} \frac{b}{b-1} \varepsilon^{-2}\right)$$

(4.14)

$$\sim \frac{b}{2}\left(\frac{b}{b-1}\right)^{1/(b-1)} \widetilde{\psi}(\varepsilon^{-2}) \qquad \text{as } \varepsilon \to 0,$$

and therefore $R$ is regularly varying at zero of index $-2/(b-1)$. Finally, by ([2.8](#)), $R(e_n) \sim a_{r(e_n)} \sim r(e_n)b/2$ as $n \to \infty$ with $a_k$ from ([4.5](#)) and thus

$$r(e_n) \sim \frac{2 \log n}{b} \qquad \text{as } n \to \infty. \qquad\qquad \square$$

PROOF OF COROLLARY 2.5. For every $n \in \mathbb{N}, c \in (0,1)$, we have

$$e_n \geq cF^{-1}\left(\log\left(\frac{n}{1-c^2}\right)\right)$$

[see Dereich, Fehringer, Matoussi and Scheutzow (2003) or Graf, Luschgy and Pagès (2003)]. Consequently,

$$F\left(\frac{e_n}{c}\right) \leq \log\left(\frac{n}{1-c^2}\right) \sim \log n \qquad \text{as } n \to \infty.$$

By Corollary 2.4, this implies

$$F\left(\frac{\varepsilon}{c}\right) \leq F\left(\frac{e_{N(\varepsilon)}}{c}\right) \lesssim \log N(\varepsilon) \sim R(\varepsilon)$$

and thus

$$F(\varepsilon) \lesssim R(c\varepsilon) \sim c^{-2/(b-1)}R(\varepsilon) \qquad \text{as } \varepsilon \to 0.$$

Letting $c \to 1$ yields the assertion. $\quad\square$

**Acknowledgment.** We wish to thank Bianca Krämer for doing some computations.

## REFERENCES

BERGER, T. (1971). *Rate Distortion Theory.* Prentice-Hall, Englewood Cliffs, NJ. MR408988

BINGHAM, N. H., GOLDIE, C. M. and TEUGELS, J. L. (1987). *Regular Variation.* Cambridge Univ. Press. MR898871

BINIA, J. (1974). On the $\varepsilon$-entropy of certain Gaussian processes. *IEEE Trans. Inform. Theory* **20** 190–196. MR384319

BRONSKI, J. C. (2003). Small ball constants and tight eigenvalue asymptotics for fractional Brownian motions. *J. Theoret. Probab.* **16** 87–100. MR1956822

DEMBO, A. and ZEITOUNI, O. (1998). *Large Deviations Techniques and Applications*, 2nd ed. Springer, New York. MR1619036

DEREICH, S. (2003). Small ball probabilities around random centers of Gaussian measures and applications to quantization. *J. Theoret. Probab.* **16** 427–449. MR1982037

DEREICH, S., FEHRINGER, F., MATOUSSI, A. and SCHEUTZOW, M. (2003). On the link between small ball probabilities and the quantization problem for Gaussian measures on Banach spaces. *J. Theoret. Probab.* **16** 249–265. MR1956830

DONOHO, D. L. (2000). Counting bits with Kolmogorov and Shannon. Technical Report 38, Stanford Univ.

FREEDMAN, D. (1999). On the Bernstein-von Mises theorem with infinite–dimensional parameters. *Ann. Statist.* **27** 1119–1140. MR1740119

GAO, F., HANNING, J. and TORCASO, F. (2003). Integrated Brownian motions and exact $L_2$-small balls. *Ann. Probab.* **31** 1320–1337. MR1989435

GERSHO, A. and GRAY, R. M. (1992). *Vector Quantization and Signal Compression.* Kluwer, Boston.

GRAF, S. and LUSCHGY, H. (2000). *Foundations of Quantization for Probability Distributions. Lecture Notes in Math.* **1730**. Springer, Berlin. MR1764176

GRAF, S., LUSCHGY, H. and PAGÈS, G. (2003). Functional quantization and small ball probabilities for Gaussian processes. *J. Theoret. Probab.* **16** 1047–1062. MR2033197

GRAY, R. M. and NEUHOFF, D. L. (1998). Quantization. *IEEE Trans. Inform. Theory* **44** 2325–2383. MR1658787

IHARA, S. (1970). On $\varepsilon$-entropy of equivalent Gaussian processes. *Nagoya Math. J.* **37** 121–130. MR258107

IHARA, S. (1993). *Information Theory*. World Scientific, Singapore. MR1249933

KOLMOGOROV, A. N. (1956). On the Shannon theory of information transmission in the case of continuous signals. *IRE Trans. Inform. Theory* **2** 102–108.

LUSCHGY, H. and PAGÈS, G. (2002). Functional quantization of Gaussian processes. *J. Funct. Anal.* **196** 486–531. MR1943099

PAPAGEORGIOU, A. and WASILKOWSKI, G. W. (1990). On the average complexity of multivariate problems. *J. Complexity* **6** 1–23. MR1048027

RITTER, K. (2000). *Average-Case Anslysis of Numerical Problems. Lecture Notes in Math.* **1733**. Springer, Berlin. MR1763973

ROSENBLATT, M. (1963). Some results on the asymptotic behavior of eigenvalues for a class of integral equations with translation kernel. *J. Math. Mech.* **12** 619–628. MR150551

ROSENFELD, A. and KAK, A. C. (1976). *Digital Picture Processing*. Academic Press, New York.

SHANNON, C. E. and WEAVER, W. (1949). *The Mathematical Theory of Communication*. Univ. Illinois Press, Urbana. MR32134

STEIN, M. L. (1999). *Interpolation of Spatial Data*. Springer, New York. MR1697409

FB IV-MATHEMATIK
UNIVERSITÄT TRIER
D-54286 TRIER
BR DEUTSCHLAND
E-MAIL: luschgy@uni-trier.de

LABORATOIRE DE PROBABILITÉS
ET MODÈLES ALÉATOIRES
UNIVERSITÉ PARIS 6
CNRS-UMR 7599
4 PLACE JUSSIEU CASE 188
F-75252 PARIS CEDEX 05
FRANCE
E-MAIL: gpa@ccr.jussieu.fr